\documentclass{amsart}
\usepackage{mathrsfs}
\usepackage{amsfonts}
\usepackage[all]{xy}
\usepackage{latexsym,amsfonts,amssymb,amsmath,amsthm,amscd,amsmath,graphicx}

\usepackage[all]{xy}

\usepackage{ifthen}

\usepackage[T1]{fontenc}
\usepackage{subsupscripts}

\newboolean{workmode}
\setboolean{workmode}{false}

\newboolean{sections}
\setboolean{sections}{true}

\input xy
\xyoption{all}

\makeatletter
\def\th@plain{%
  \thm@notefont{}
  \itshape 
}
\def\th@definition{%
  \thm@notefont{}
  \normalfont 
}
\makeatother

\newcommand{\hypref}[2]{{\hyperref[#1]{#2~\ref{#1}}}}

\newcommand{\ifwork}[1]{\ifthenelse{\boolean{workmode}}{#1}{}}
\newcommand{\comment}[1]{}
\newcommand{\mute}[1]{}
\newcommand{\printname}[1]{}

\ifwork{
\renewcommand{\comment}[1]{{\marginpar{*}\ \scriptsize{#1}\ }}
\renewcommand{\mute}[1]{{\scriptsize \ #1\ }\marginpar{\scriptsize muted}}
\renewcommand{\printname}[1]
    {\smash{\makebox[0pt]{\hspace{-1.0in}\raisebox{8pt}{\tiny #1}}}}
}

\newcommand{\ifsection}[2]{\ifthenelse{\boolean{sections}}{#1}{#2}}

\theoremstyle{plain}
\ifsection{
    \newtheorem{theorem}{Theorem}[section]
}
{
    
}

\theoremstyle{definition}

\newtheorem{thm}{Theorem}[section]
\newtheorem{cor}[thm]{Corollary}
\newtheorem{lem}[thm]{Lemma}

\theoremstyle{definition}
\newtheorem{defn}[thm]{Definition}

\theoremstyle{remark}
\newtheorem{rem}[thm]{Remark}
\numberwithin{equation}{section}

\newcommand{\sA}{{\mathcal A}}

\newcommand{\sE}{{\mathcal E}}
\newcommand{\sF}{{\mathcal F}}

\newcommand{\sK}{{\mathcal K}}

\newcommand{\sM}{{\mathcal M}}

\newcommand{\sO}{{\mathcal O}}

\newcommand{\sV}{{\mathcal V}}

\newcommand{\G}{{\mathbb G}}

\newcommand{\N}{{\mathbb N}}
\renewcommand{\P}{{\mathbb P}}
\newcommand{\Q}{{\mathbb Q}}
\newcommand{\Z}{{\mathbb Z}}

\newcommand{\x}{\xrightarrow}

\numberwithin{equation}{section}

\begin {document}
\topmargin= -.2in \baselineskip=20pt

\title{Special $L$-values of abelian $t$-modules}
\author{Jiangxue Fang}
\address{Department of Mathematics,
Capital Normal University, Beijing 100148, P.R. China} \email{fangjiangxue@gmail.com}
\date{}

\maketitle
\begin{abstract}
We prove a formula for the $\infty$-adic special $L$-value of abelian $t$-modules. This gives function field analogues of the
class number formula. We also express it in terms of the extension groups of shtukas.
\end{abstract}

\section{{Introduction and statement of the main results}}
\noindent Let $k$ be a finite field of $q$ elements. Without mention to the contrary schemes are understood to be over $k$ and tensor products over $k$.
For any finite $k[t]$-module $M$, let
$$|M|={\rm det}_{k[t]}\Big(1\otimes t-t\otimes 1,\;M\otimes_kk[t]\Big),$$
where the $k[t]$-module structure on $M\otimes_kk[t]$ is given by $k[t]$.
Let $R$ be the integral closure of $k[t]$ in a finite extension $K$ of $k(t)$ and let $z$ be the image of $t$ in $R$.
For any $n\geq1$, consider in $k((t^{-1}))$, the infinite sum
$$\sum_{I}\frac{1}{|R/I|^n}$$ where $I$ runs over all nonzero ideals of $R$,
converges to an element of $1+t^{-1}k[[t^{-1}]]$, which is denoted by $\zeta(R,\,n)$.
By the product formula,
$$\zeta(R,\,n)=\prod_{\mathfrak m\in{\rm Max}(R)}\Big(1-\frac{1}{|R/\mathfrak m|^n}\Big)^{-1}$$
where ${\rm Max}(R)$ is the set of all maximal ideals of $R$.
\begin{defn}
The $n$-th tensor power of Carlitz module is the functor
$$C^{\otimes n}:\{R\hbox{-algebras}\}\to\{k[t]\hbox{-modules}\}$$
that associates each $R$-algebra $B$ to a $k[t]$-module $C^{\otimes n}(B)$ whose underlying $k$-vector space is $B^n$ and whose
$k[t]$-module structure is given by
$$\varphi:k[t]\to{\rm End}_k(B^n),\;\;\varphi(t)(x_1,\ldots,x_{n-1},x_n)=(zx_1+x_2,\ldots,zx_{n-1}+x_n,zx_n+x_1^q)$$
for any $x_1,\ldots,x_n\in B.$
\end{defn}
\begin{lem}
We have
$$|C^{\otimes n}(R/\mathfrak m)|=|R/\mathfrak m|^n-1,\;$$
and then
$$\zeta(R,\,n)=\prod_{\mathfrak m\in{\rm Max}(R)}\frac{|R/\mathfrak m|^n}{|C^{\otimes n}(R/\mathfrak m)|}.$$
\end{lem}

For any matrix $(a_{ij})$ over a $k$-algebra, denote by $(a_{ij})^{(q^s)}=(a_{ij}^{q^s})$. Let $M_n(R)$ be the ring of
$n\times n$-matrixes over $R$ and let $M_n(R)\{\tau\}$ be the ring over $M_n(R)$ generated by $\tau$ with the relation $\tau P=P^{(q)}\tau$ for any $P\in M_n(R)$.

\begin{defn}
An abelian $t$-module over $R$ is a $k[t]$-module scheme $E$ over $R$ whose underlying $k$-vector space scheme is
isomorphic to $\G_{\rm a}^n$ for some $n$ and the $k[t]$-module
structure on $E$ is given by
$$\varphi_E:k[t]\to{\rm End}_k(E)=M_n(R)\{\tau\},\;t\mapsto\sum_{s=0}^rA_s\tau^s$$
for some $A_0,\ldots,A_r\in M_n(R)$ such that $(A_0-zI_n)^n=0$. The integer $n$ is called the dimension of $E$.
\end{defn}

For an abelian $t$-module $E$ over $R$, define the $\infty$-adic special $L$-value of $E$ over $R$ by
$$L(E/R)=\prod_{\mathfrak m\in{\rm Max}(R)}\frac{|{\rm Lie}(E)(R/\mathfrak m)|}{|E(R/\mathfrak m)|}\in 1+t^{-1}k[[t^{-1}]].$$
The $n$-th tensor power of Carlitz module is an abelian $t$-module and $\zeta(R,\,n)=L(C^{\otimes n}/R)$.

Let $K_\infty=K\otimes_{k(t)}k((t^{-1}))$. There exists a unique power series
$$\exp_EX=\sum_{s\geq0}e_sX^{(q^s)}$$
with $X=(X_1,\ldots,X_n)^T$ and $e_s\in M_n(K_\infty)$ such that $e_0=I_n$ and
$$\exp_E(A_0X)=\sum_{s=0}^rA_s(\exp_EX)^{(q^s)}.$$
The valuation of $k((t^{-1}))$ extends to a valuation ${val}$ on $K_\infty$, $K_\infty^n$ and $M_n(K_\infty)$ in the natural way. By Proposition 2.14 of \cite{A}, we have
$$\lim_{s\to\infty}\frac{{val}(e_s)}{q^s}=+\infty$$
and we get a continuous and open $k[t]$-linear map
$$\exp_E:{\rm Lie}(E)(K_\infty)\to E(K_\infty).$$

\begin{lem}\label{fu}
We have $A_0^{q^n}=z^{q^n}I_n$ and
$\inf_{s\in\Z}val(A_0^s)+s>-\infty.$
\end{lem}
\begin{proof}
Let $P=A_0-zI_n.$ Then $P^n=0$. We have $A_0^{q^n}=(zI_n+P)^{q^n}=z^{q^n}I_n$. For any $s\in\Z$, there exists $t\in\Z$ such that $0\leq s-tq^n<q^n$.
Then $A_0^s=A_0^{s-tq^n}A_0^{tq^n}=z^{tq^n}A_0^{s-tq^n}$ and hence
$$val(A_0^s)+s=val(z^{tq^n}A_0^{s-tq^n})+s=val(A_0^{s-tq^n})+s-tq^n.$$
This proves $\inf_{s\in\Z}val(A_0^s)+s>-\infty$.
\end{proof}
Then $\sum_{s\ll+\infty}a_sA_0^s(x)$ converges for any $\sum_{s\ll+\infty}a_st^s\in k((t^{-1}))$ with $a_s\in k$ and any $x\in {\rm Lie}(E)(K_\infty)$.
This gives a $k((t^{-1}))$-vector space structure on ${\rm Lie}(E)(K_\infty)$.

\begin{defn}
Suppose $S$ is a commutative ring. A perfect complex of $S$-modules is a bounded complex of projective $S$-modules of finite type.
 Let $D^{\rm per}(S)$ be the full subcategory of $D(S)$ consisting of all objects which can be
 represented by a perfect complex. For any perfect complex $C=(C^i)_{i\in\Z}$, the determinant $\det(C)$ of $C$ is defined by $\bigotimes_{i\in\Z}\det(C^i)^{(-1)^i}$.
For any two isomorphisms $f,\,g:\det(C_1)\simeq \det(C_2)$ for some $C_1,\,C_2\in D^{\rm perf}(S)$, denote by $[f:g]$
the element of $S^\times$ such that $f=[f:g]g$.
\end{defn}

\begin{defn}\label{16} Let $V$ be a finite dimensional $k((t^{-1}))$-vector space. A lattice in $V$ is a finite free sub-$k[t]$-module $\Lambda$ of $V$ such that
the natural morphism $\Lambda\otimes_{k[t]}k((t^{-1}))\to V$ is isomorphic.
Let $\Lambda_1$ and $\Lambda_2$ be two lattices of $V$.
Any isomorphism $\det(\Lambda_1)\simeq\det(\Lambda_2)$ defines an isomorphism $\det\Big(\Lambda_1\otimes_{k[t]}k((t^{-1}))\Big)\simeq
\det\Big(\Lambda_2\otimes_{k[t]}k((t^{-1}))\Big).$ Let $\rho$ be the composition
$$\det(V)\simeq\det\Big(\Lambda_1\otimes_{k[t]}k((t^{-1}))\Big)\simeq
\det\Big(\Lambda_2\otimes_{k[t]}k((t^{-1}))\Big)\simeq\det(V).$$
Then the image of $[\rho:1]$ in $k((t^{-1}))^\times/k^\times$ does not depend on the choice of the isomorphism $\det(\Lambda_1)\simeq\det(\Lambda_2)$.
Denote by $[\Lambda_1:\Lambda_2]$ the monic representative
of $[\rho:1]$ in $k((t^{-1}))^\times$.
\end{defn}

Define
$$H(E/R)=\frac{E(K_\infty)}{\exp_E({\rm Lie}(E)(K_\infty))+E(R)}.$$
In this paper, we prove the following theorem.
\begin{thm}\label{18}
Let $E$ be an abelian $t$-module over $R$. Then ${\rm Lie}(E)(R)$ and $\exp_E^{-1}(E(R))$ are lattices of ${\rm Lie}(E)(K_\infty)$ and $H(E/R)$
is a finite $k[t]$-module. Moreover,
$$L(E/R)=[{\rm Lie}(E)(R):\exp_E^{-1}(E(R))]\cdot|H(E/R)|\in k((t^{-1}))^\times.$$
\end{thm}

\begin{rem}
The statement of Theorem \ref{18} is similar to the statement of the class number formula. Let $F$ be a number field and $\sO_F$ be the ring of integers in $F$.
Consider the exponential map
$$\exp:(\sO_F/\Z)\otimes_\mathbb Z\mathbb R\to(\sO_F\otimes_\mathbb Z\mathbb R)^\times/\mathbb R^\times_{>0}.$$
Both $\sO_F/\mathbb Z$ and $\exp^{-1}(\sO_F^\times)$ are lattices of $(\sO_F/\mathbb Z)\otimes_\mathbb Z\mathbb R$. The class number formula is printed as
\begin{eqnarray*}
&&\lim_{s\to1}(s-1)\zeta_F(s)=\frac{2^{r_1}(2\pi)^{r_2}R_Fh_F}{w_F\sqrt{|\Delta_F|}}=\lambda\cdot[\sO_F/\mathbb\Z:\exp^{-1}(\sO_F^\times)]
\end{eqnarray*}
for some $\lambda\in\Q^\times.$
\end{rem}

\begin{rem}
This theorem generalizes results of Taelman \cite{T} for Drinfeld modules which are abelian $t$-modules of dimension one and results of Anderson and Thakur \cite{AT} for $\zeta(k[t],\,n).$ In \cite{T}, Taelman use a variant of Anderson's trace formula and we use the determinants of extension group of shtukas considered in \cite{L}.
\end{rem}

\begin{defn}
Let $T$ and $S$ be two $k$-schemes. Denote by $F_S:S\to S$ the morphism defined by $\sO_S\to\sO_S,\;a\mapsto a^q$.  A $T$-shtuka on $S$ is a diagram
$$M\x{i}M'\xleftarrow{j}M$$
where $M$ and $M'$ are quasi-coherent $\sO_{S\times T}$-modules, $i$ is $\sO_{S\times T}$-linear and $j$ is $F_S\times 1$-linear. The category of $T$-shtukas on $S$ is abelian. We can define the extension groups ${\rm Ext}^\bullet(M_1,\;M_2)$ of two $T$-shtukas $M_1$ and $M_2$ on $S$.
The unit $T$-shtuka $\textbf{1}_{S\times T}$ on $S$ is defined to be
$$\textbf{1}_{S\times T}=[\sO_{S\times T}\x{1}\sO_{S\times T}\xleftarrow{F_S\times{1}}\sO_{S\times T}].$$
\end{defn}

The extension $k(t)\hookrightarrow K$ defines a surjective morphism $X\to\P^1$ of smooth projective curves. Let $Y={\rm Spec}\; R$.
For any integer $d$, denote by $\sO_X(d\infty)$ the pullback of
$\sO_{\P^1}(d\infty)$ via $X\to\P^1$.
For $d\gg0$, $A_1,\ldots, A_r\in M_n(\Gamma(X,\,\sO_X((q-1)d\infty)))$. Thus $A_1,\ldots,A_r$
define morphisms $\sO_X(-qd\infty)^n={F}_X^*\sO_X(-d\infty)^n\to\sO_X(-d\infty)^n$ of vector bundles on $X$. The morphism $\sO_X\to\sO_X,\;a\mapsto a^q$
defines a ${F}_X$-linear map $\tau:\sO_X(-d\infty)^n\to\sO_X(-qd\infty)^n$. Then we get $F_X$-linear maps $A_1\tau,\ldots,A_r\tau:\sO_X(-d\infty)^n\to\sO_X(-d\infty)^n$ for $d\gg0$.

For any quasi-coherent $\sO_X$-module $\sM$ on $X$, denote by $\sM[t]$ the inverse image of $\sM$ under the projection $X\otimes k[t]\to X$.

\begin{defn}\label{ke}
Let $d\gg0$ and $e\geq1$ such that $A_0\in M_n(\Gamma(X,\,\sO_X(e\infty))).$

(1) Let $\widetilde\sE$ be the ${\rm Spec}\,k[t]$-shtuka on $X$:
\begin{eqnarray*}
\Big(\sO_X(-d\infty)^n[t]\Big)^r\x{\widetilde i}\Big(\sO_X(-d\infty)^n[t]\Big)^{r-1}\bigoplus \sO_X((e-d)\infty)^n[t]
\xleftarrow{\widetilde j}\Big(\sO_X(-d\infty)^n[t]\Big)^r
\end{eqnarray*}
where
\begin{eqnarray*}
&&\widetilde i(x_1,\ldots,x_{r-1},x_r)=(-x_2,\ldots,-x_r,(t-A_0)(x_1));\\
&&\widetilde j(x_1,\ldots,x_{r-1},x_r)=(-\tau(x_1),\ldots,-\tau(x_{r-1}),\sum_{s=1}^rA_s\tau(x_s))
\end{eqnarray*}
for any $x_1,\ldots,x_r\in\sO_X(-d\infty)^n[t]$.

(2) Let
\begin{eqnarray*}
\sE^\bullet=\Big(\sO_X(-d\infty)^n[t]\x{t-\sum_{s=0}^rA_s\tau^s}\sO_X((e-d)\infty)^n[t]\Big).
\end{eqnarray*}

(3) For any place $w$ of $X$, let $K_w$ be the completion of $K$ at the place $w$. Define a sheaf of $k[t]$-modules $\mathfrak E$ on $X$ by
\begin{eqnarray*}
\mathfrak E(U)=\Big\{(x,\,(\gamma_w)_w)\in E(\sO_X(U))\times\prod_{w\in U\backslash Y}{\rm Lie}(E)(K_w)\;|\,\exp_E(\gamma_w)=x\hbox{ for any }w\Big\}
\end{eqnarray*}
for any open subset $U$ of $X$.
\end{defn}

\begin{thm}\label{yang}
Suppose $d\gg0$.

(1) The complex $\sE^\bullet[1]$ is a sheaf and there is a natural surjective morphism of sheaves $\sE^\bullet[1]\to\mathfrak E$ on $X$ whose kernel is a skyscraper sheaf $\sK$ on $X$ supported on $X-Y$. We have
\begin{eqnarray*}
&&H^0(X,\,\mathfrak E)= \exp_E^{-1}(E(R));\\
&&H^1(X,\,\mathfrak E)= H(E/R).
\end{eqnarray*}

(2) We have
\begin{eqnarray*}
&&{\rm Ext}(\textbf{1}_{X\otimes k[t]},\;\widetilde\sE)\simeq R\Gamma(X,\,\sE^\bullet);\\
&&{\rm Ext}^2(\textbf{1}_{X\otimes k[t]},\;\widetilde\sE)\simeq H(E/R);\\
&&{\rm Ext}^s(\textbf{1}_{X\otimes k[t]},\;\widetilde\sE)=0\hbox{ for }
s\neq 1\hbox{ and }2.
\end{eqnarray*}
The surjective morphism $\sE^\bullet[1]\to\mathfrak E$ induces a surjective homomorphism ${\rm Ext}^1(\textbf{1}_{X\otimes k[t]},\;\widetilde\sE)\to \exp_E^{-1}(E(R))$ whose kernel is $\sK$.

(3) Suppose furthermore $A_0\in M_n(\Gamma(X,\,\sO_X(\infty)))$. One can take $e=1$. Then
$\sE^\bullet[1]\simeq\mathfrak E$ and
\begin{eqnarray*}
&&{\rm Ext}^1(\textbf{1}_{X\otimes k[t]},\;\widetilde\sE)\simeq\exp_E^{-1}(E(R)).
\end{eqnarray*}
\end{thm}

The paper is organized as follows. In section 2, we express the $v$-adic special values of shtukas in terms of the determinant of the extension groups of shtukas under some local
analytic conditions. In section 3, we prove the class number formula for an abelian $t$-module. In section 4, we express $\exp_E^{-1}(E(R))$ and $H(E/R)$
in terms of the extension groups of some shtuka.

\textbf{Acknowledgements.} I would like to thank Lenny Taelman for his excellent course about the Woods Hole trace formula in Morningside Center. 
My research is supported by the NSFC.

\section{$v$-adic $L$-values of shtukas on curves}
\noindent In this section, let $T={\rm Spec}\;A$ be a smooth affine curve over $k$. Let $i:\sE\to\sE'$ be a morphism of vector bundles on $X\times T$ such that $i$ is isomorphic at the generic point. Let $j_0:\sE\to\sE$ and $j_1,\,\ldots,j_r:\sE\to\sE'$ be $F_X\times 1$-linear maps of $\sO_{X\times T}$-modules and let $j=\sum_{s=1}^rj_sj_0^{s-1}:\sE\to\sE'$.

Let $Z(\det(i))$ be the zeros of $\det(i)$ in $X\times T$. Fix $v\in|T|$ such that $Z(\det(i))\cap X\times\{v\}$ is finite. Take a finite subset $S$ of $|X|$ such that $S\times\{v\}\supset Z(\det(i))\cap X\times\{v\}$. Let $A_v$ be the completion of $A$ at $v$ and choose a uniformizer of $A_v$ which is also denoted by $v$.
For any $x\in|X|-S$, let $i_x,\,j_x,\,(j_1)_x,\ldots,(j_r)_x:\sE_x\to\sE'_x$ and $(j_0)_x:\sE_x\to\sE_x$ be the restriction of $i,\,j,\,j_1,\ldots,j_r$ and $j_0$
on ${\rm Spec}\;k(x)\otimes A_v$, respectively.
Then $i_x:\sE_x\to\sE_x'$ is isomorphic for any $x\in|X|-S$.

\begin{lem}
Define the $v$-adic $L$-function away from $S$ of the diagram
$\sE\x{i}\sE'\xleftarrow{j}\sE$
to be
$$L_{v}(X-S,\,(\sE,\sE',i,j),T)=\prod_{x\in|X|-S}{\rm det}_{A_v}\Big(1-\sum_{s=1}^rT^si_x^{-1}(j_s)_x(j_0)_x^{s-1},\,\sE_x\Big)^{-1}\in 1+TA_v[[T]].$$
Then $L_{v}(X-S,\,(\sE,\sE',i,j),T)\in 1+TA_v\langle\langle T\rangle\rangle$.
\end{lem}
\begin{proof}
Consider the $T$-shtuka
$$\sE^r\x{\widetilde i}\sE^{r-1}\oplus\sE'\xleftarrow{\widetilde j}\sE^r$$
on $X$, where
\begin{eqnarray*}
&& \widetilde i(x_1,\ldots,x_{r-1}, x_r)=(-x_2,\ldots,-x_r,i(x_1));\\
&& \widetilde j(x_1,\ldots,x_{r-1},x_r)=(-j_0(x_1),\ldots,-j_0(x_{r-1}),\sum_{s=1}^rj_s(x_s))
\end{eqnarray*}
for any $x_1,\ldots,x_r\in\sE$.
For any $x\in|X|-S$, let $\widetilde i_x$ and $\widetilde j_x:\sE_x^r\to\sE_x^{r-1}\oplus\sE_x'$ be the restriction of $\widetilde i$ and $\widetilde j$
on ${\rm Spec}\;k(x)\otimes A_v$, respectively.
We have $Z(\det(\widetilde i))=Z(\det(i))$ and then $\widetilde i_x:\sE_x^r\to\sE^{r-1}_x\oplus\sE_x'$ is isomorphic for any $x\in|X|-S$.
The $v$-adic $L$-function away from $S$ of the $T$-shtuka $(\sE^r,\sE^{r-1}\oplus\sE',\widetilde i,\widetilde j)$ on $X$ is defined to be
$$L_{v}(X-S,\,(\sE^r,\sE^{r-1}\oplus\sE',\widetilde i,\widetilde j),\,T)=\prod_{x\in|X|-S}{\rm det}_{A_v}(1-T(\widetilde i_x)^{-1}\widetilde j_x,\,\sE_x^r)^{-1}\in 1+TA_v[[T]].$$
By \cite{L} Proposition 3.2 (a), $L_{v}(X-S,\,(\sE^r,\sE^{r-1}\oplus\sE',\widetilde i,\widetilde j),\,T)\in 1+TA_v\langle\langle T\rangle\rangle$.
By Lemma \ref{d} (3), we have
\begin{eqnarray}\label{xiu}
L_{v}(X-S,\,(\sE,\sE',i,j),T)=L_{v}(X-S,\,(\sE^r,\sE^{r-1}\oplus\sE',\widetilde i,\widetilde j),\,T)\in 1+TA_v\langle\langle T\rangle\rangle.
\end{eqnarray}
\end{proof}

\begin{defn}
We have
$$L_{v}(X-S,\,(\sE,\sE',i,j),T)=(1-T)^sL_{v}^*(X-S,\,(\sE,\sE',i,j),T)$$ for some $s\in\N$ and $L_{v}^*(X-S,\,(\sE,\sE',i,j),T)\in 1+TA_v\langle\langle T\rangle\rangle$ such that
$L_{v}^*(X-S,\,(\sE,\sE',i,j),1)\in A_v^\times$. The $v$-adic special $L$-value $L_{v}^*(X-S,\,(\sE,\sE',i,j))$
of the diagram $(\sE,\sE',i,j)$ away from $S$ is defined to be $L_{v}^*(X-S,(\sE,\sE',i,j),1)\in A_v^\times.$
\end{defn}
\begin{lem}\label{d}
Let $\sA$ be an abelian category. Let $M,\,M'\in{\rm Ob}\sA$ and let $i,j_1,\ldots, j_r\in{\rm Hom}(M,M')$ and $j_0\in{\rm End}(M)$.
Define $\widetilde i,\;\widetilde j:M^r\to M^{r-1}\oplus M'$, $\phi_1,\,\sigma:M^r\to M^r$ and $\phi_2:M^{r-1}\oplus M'\to M^{r-1}\oplus M'$ by
\begin{eqnarray*}
&&\widetilde i(x_1,\ldots,x_{r-1}, x_r)=(-x_2,\ldots,-x_r,i(x_1));\\
&&\widetilde j(x_1,\ldots,x_{r-1},x_r)=(-j_0(x_1),\ldots,-j_0(x_{r-1}),\sum_{s=1}^rj_s(x_s))\\
&&\phi_1(x_1,\ldots,x_{r-1},x_r)=(j_0(x_1)-x_2,\ldots,j_0(x_{r-1})-x_r,x_1);\\
&&\sigma(x_1,\ldots,x_{r-1},x_r)=(-x_2,\ldots,-x_r,x_1);
\\&&\phi_2(x_1,\ldots,x_{r-1},x')=(x_1,\ldots,x_{r-1},x'-\sum_{s+t\leq r,\;1\leq s,t}j_{s+t}j_0^{s-1}(x_t))
\end{eqnarray*}
for any $x_1,\ldots,x_r\in M$ and $x'\in M'$. Let $p_r:M^r\to M$ and $p_r:M^{r-1}\oplus M'\to M'$ be the projections to the last factors.

(1) The commutative diagram
\[\xymatrix{M^r\ar[d]^{\widetilde i}\ar[r]^{\sigma}&M^r\ar[r]^{p_r}\ar[d]^{{\rm id}_{M^{r-1}}\oplus i}&M\ar[d]^i\\
M^{r-1}\oplus M'\ar[r]^{\rm id}&M^{r-1}\oplus M'\ar[r]^{p_r}&M'}\]
defines quasi-isomorphisms
$$\Big(M^r\x{\widetilde i}M^{r-1}\oplus M'\Big)\simeq\Big(M^r\x{{\rm id}_{M^{r-1}}\oplus i}M^{r-1}\oplus M'\Big)\simeq\Big(M\x{i}M'\Big).$$

(2) Let $j=\sum_{s=1}^rj_sj_0^{s-1}$.
The commutative diagram
\[\xymatrix{M^r\ar[d]^{\widetilde i-\widetilde j}\ar[r]^{\phi_1}&M^r\ar[r]^{p_r}\ar[d]^{{\rm id}_{M^{r-1}}\oplus i-j}&M\ar[d]^{i-j}\\
M^{r-1}\oplus M'\ar[r]^{\phi_2}&M^{r-1}\oplus M'\ar[r]^{p_r}&M'}\]
defines quasi-isomorphisms
$$\Big(M^r\x{\widetilde i-\widetilde j} M^{r-1}\oplus M'\Big)\simeq\Big(M^r\x{{\rm id}_{M^{r-1}}\oplus i-j}M^{r-1}\oplus M'\Big)\simeq\Big(M\x{i-j}M'\Big).$$

(3) Suppose $\sA$ is the category of $S$-modules for some commutative ring $S$. Suppose $M$ and $M'$ are finite free $S$-modules and $i:M\to M'$ is an isomorphism.
Then
\begin{eqnarray*}
&&{\rm det}_S(1-\widetilde i^{-1}\widetilde j,\,M^r)={\rm det}_S(1-i^{-1}j,\,M);\\
&&{\rm det}_S(1-T\widetilde i^{-1}\widetilde j,\,M^r)={\rm det}_S(1-\sum_{s=1}^rT^si^{-1}j_s(j_0)^{s-1},\,M).
\end{eqnarray*} \end{lem}
\begin{proof}
(1) is a special case of (2) when $j_0,\ldots, j_r$ are zero map. Define $\phi:M^r\to M^{r-1}\oplus M'$ by
$$\phi(x_1,\ldots,x_{r-1},x_r)=(x_1,\ldots,x_{r-1},i(x_r)-\sum_{s=1}^rj_sj_0^{s-1}(x_r)+\sum_{s+t\leq r,\;1\leq s,t}j_{s+t}j_0^{s-1}(x_t)).$$
Then (2) follows from the fact $\phi\phi_1=\widetilde i-\widetilde j,\;\phi_2\phi={\rm id}\oplus i-j$ and the bijectivity of $\phi_1$ and $\phi_2$.
The first equality of (3) follows form the fact $\det(\sigma)=\det(\phi_1)=\det(\phi_2)=1$. Applying this results to $k[T]$-linear maps
$i\otimes1:M\otimes_k k[T]\to M'\otimes_kk[T]$, $Tj_0:M\otimes_kk[T]\to M\otimes_kk[T]$ and $Tj_1,\ldots,Tj_r:M\otimes_kk[T]\to M'\otimes_kk[T]$, we get the second equality of (3).
\end{proof}

Suppose $X-S={\rm Spec}\;R_S$. Let $\sM$ and $\sM'$ be the $R_S\widehat\otimes A_v$-modules defined by $\sE$ and $\sE'$. For any $w\in S$, let $\sO_w$ be the ring of integers in $K_w$. Choose a uniformizer of $\sO_w$ and denote it also by $w$. Let $\sM_w$ and $\sM_w'$ be the $\sO_w\widehat\otimes A_v$-modules defined by $\sE$ and $\sE'$. Let $\sV={\rm Hom}_{R_S\widehat\otimes A_v}(\sM,\Omega_{R_S}\widehat\otimes A_v)$ and $\sV'={\rm Hom}_{R_S\widehat\otimes A_v}(\sM',\Omega_{R_S}\widehat\otimes A_v)$.
We have the commutative diagram of $A_v$-modules
\[\xymatrix{R\Gamma(X,\,\sE)\otimes_AA_v\ar[r]\ar@<0.5ex>[d]_{i-j~~}\ar@<-0.5ex>[d]
^{~~i}&\bigoplus_{w\in S}\sM_w\ar[r]\ar@<0.5ex>[d]_{i-j~~}\ar@<-0.5ex>[d]^{~~i}&{\rm Hom}_{A_v}(\sV,\,A_v)\ar@<0.5ex>[d]_{i-j~~}\ar@<-0.5ex>[d]^{~~i}\\
R\Gamma(X,\,\sE')\otimes_AA_v\ar[r]&\bigoplus_{w\in S}\sM_w'\ar[r]&{\rm Hom}_{A_v}(\sV',\,A_v),}\]
where the morphism $\bigoplus_{w\in S}\sM_w\to{\rm Hom}_{A_v}(\sV,\,A_v)$ (resp. $\bigoplus_{w\in S}\sM'_w\to{\rm Hom}_{A_v}(\sV',\,A_v)$)
associates each $(f_w)\in\bigoplus\sM_w$ and $g\in\sV$ (resp. $(f_w)\in\bigoplus\sM_w'$ and $g\in\sV'$) to the sum of residue of $\langle g,\,f_w\rangle$ at $w$.

Since $S\times\{v\}\supset Z(\det(i))\cap X\times\{v\}$, then $i:{\rm Hom}_{A_v}(\sV,\,A_v)\simeq{\rm Hom}_{A_v}(\sV',\,A_v)$ is an isomorphism. It defines a quasi-isomorphism
$$R\Gamma(X,\,\sE\x{i}\sE')\otimes_AA_v\simeq\bigoplus_{w\in S}(\sM_w\x{i}\sM_w')$$
and an isomorphism
$$\alpha:\det\Big(R\Gamma(X,\,\sE\x{i}\sE')\otimes_AA_v\Big)\simeq\bigotimes_{w\in S}\det\Big(\sM_w\x{i}\sM_w'\Big).$$
By Lemma \ref{d} (1), $\widetilde i:{\rm Hom}_{A_v}(\sV^r,\,A_v)\simeq{\rm Hom}_{A_v}(\sV^{r-1}\oplus\sV',\,A_v)$ is an isomorphism and it defines an isomorphism
$$\widetilde\alpha:\det\Big(R\Gamma(X,\,\sE^r\x{\widetilde i}\sE^{r-1}\oplus\sE')\otimes_AA_v\Big)
\simeq\bigotimes_{w\in S}\det\Big(\sM^r_w\x{\widetilde i}\sM^{r-1}_w\oplus\sM_w'\Big).$$

Let $\varphi$ be the composition
\begin{eqnarray*}
\det\Big(R\Gamma(X,\,\sE\x{i}\sE')\otimes_AA_v\Big)&\simeq&\det\Big(R\Gamma(X,\,\sE)\otimes_AA_v\Big)
\bigotimes\det\Big(R\Gamma(X,\,\sE')\otimes_AA_v\Big)^{-1}\\
&\simeq&\det\Big(R\Gamma(X,\,\sE\x{i-j}\sE')\otimes_AA_v\Big).
\end{eqnarray*}
Similarly, we have an isomorphism
$$\widetilde\varphi:\det\Big(R\Gamma(X,\,\sE^r\x{\widetilde i}\sE^{r-1}\oplus \sE')\otimes_AA_v\Big)
\simeq\det\Big(R\Gamma(X,\,\sE^r\x{\widetilde i-\widetilde j}\sE^{r-1}\oplus\sE')\otimes_AA_v\Big).$$
By the same method in Lemma 4.3 of \cite{L}, for any $n\in\N$
$$i\hbox{ and }i-j:w^t\sM_w\to w^t\sM_w'$$
are injective and they have same image and hence same cokernel which are finite free $A_v/v^nA_v$-modules for $t$ large enough. The natural quasi-isomorphism
$$\Big(\sM_w/v^n\sM_w\x{i}\sM'_w/v^n\sM_w'\Big)\simeq\Big(\sM_w/w^t\sM_w+v^n\sM_w\x{i}(\sM_w'/v^n\sM_w')/i(w^t(\sM_w/v^n\sM_w))\Big)$$
defines an isomorphism
\begin{eqnarray}\label{a}
&&\det\Big(\sM_w/v^n\sM_w\x{i}\sM'_w/v^n\sM_w'\Big)\\\nonumber&\simeq&\det\Big(\sM_w/w^t\sM_w+v^n\sM_w\x{i}(\sM_w'/v^n\sM_w')/i(w^t(\sM_w/v^n\sM_w))\Big).
\end{eqnarray}
The natural quasi-isomorphism
\begin{eqnarray*}&&\Big(\sM_w/v^n\sM_w\x{i-j}\sM'_w/v^n\sM_w'\Big)\\
&\simeq&\Big(\sM_w/w^t\sM_w+v^n\sM_w\x{i-j}(\sM_w'/v^n\sM_w')/(i-j)(w^t(\sM_w/v^n\sM_w))\Big)
\end{eqnarray*}
defines an isomorphism
\begin{eqnarray}\label{b}
&&\det\Big(\sM_w/v^n\sM_w\x{i-j}\sM'_w/v^n\sM_w'\Big)\\\nonumber&\simeq&\det\Big(\sM_w/w^t\sM_w+v^n\sM_w\x{i-j}(\sM_w'/v^n\sM_w')/(i-j)(w^t(\sM_w/v^n\sM_w))\Big).
\end{eqnarray}
Similar as $\varphi$, we have an isomorphism
\begin{eqnarray}\label{c}
&&\det\Big(\sM_w/w^t\sM_w+v^n\sM_w\x{i}(\sM_w'/v^n\sM_w')/i(w^t(\sM_w/v^n\sM_w))\Big)\\\nonumber
&\simeq&\det\Big(\sM_w/w^t\sM_w+v^n\sM_w\x{i-j}(\sM_w'/v^n\sM_w')/(i-j)(w^t(\sM_w/v^n\sM_w))\Big).
\end{eqnarray}
Then (\ref{a}), (\ref{b}) and (\ref{c}) define an isomorphism
$$\det\Big(\sM_w/v^n\sM_w\x{i}\sM'_w/v^n\sM_w'\Big)\simeq\det\Big(\sM_w/v^n\sM_w\x{i-j}\sM'_w/v^n\sM_w'\Big).$$
Taking the inverse limit, we get an isomorphism
$$\psi:\bigotimes_{w\in S}\det\Big(\sM_w\x{i}\sM'_w\Big)\simeq\bigotimes_{w\in S}\det\Big(\sM_w\x{i-j}\sM'_w\Big).$$
Using the same method, we have an isomorphism
$$\widetilde\psi:\bigotimes_{w\in S}\det\Big(\sM^r_w\x{\widetilde i}\sM_w^{r-1}\oplus\sM'_w\Big)\simeq
\bigotimes_{w\in S}\det\Big(\sM^r_w\x{\widetilde i-\widetilde j}\sM_w^{r-1}\oplus\sM'_w\Big).$$

Suppose $i-j:{\rm Hom}_{A_v}(\sV,\,A_v)\to{\rm Hom}_{A_v}(\sV',\,A_v)$ is an isomorphism. By the above commutative diagram, it defines a quasi-isomorphism
$$R\Gamma(X,\,\sE\x{i-j}\sE')\otimes_AA_v\simeq\bigoplus_{w\in S}(\sM_w\x{i-j}\sM_w')$$
and an isomorphism
$$\beta:\det\Big(R\Gamma(X,\,\sE\x{i-j}\sE')\otimes_AA_v\Big)\simeq\bigotimes_{w\in S}\det\Big(\sM_w\x{i-j}\sM_w'\Big).$$
By Lemma \ref{d} (2), $\widetilde i-\widetilde j:{\rm Hom}_{A_v}(\sV^r,\,A_v)\to{\rm Hom}_{A_v}(\sV^{r-1}\oplus\sV',\,A_v)$ is also an isomorphism, hence we get an isomorphism
$$\widetilde\beta:\det\Big(R\Gamma(X,\,\sE^r\x{\widetilde i-\widetilde j}\sE^{r-1}\oplus\sE')\otimes_AA_v\Big)\simeq
\bigotimes_{w\in S}\det\Big(\sM_w^r\x{\widetilde i-\widetilde j}\sM_w^{r-1}\oplus\sM_w'\Big).$$
\begin{lem}  \label{e}
Suppose $i-j:{\rm Hom}_{A_v}(\sV,\,A_v)\to{\rm Hom}_{A_v}(\sV',\,A_v)$ is an isomorphism. We have
$$L_{v}^*(X-S,\,(\sE,\sE',i,j))=[\psi\alpha:\beta\varphi]\in A_v^\times.$$
\end{lem}
\begin{proof}

By Lemma \ref{d}, (1) and (2), we get isomorphisms
\begin{eqnarray*}
&&\eta_i:\det\Big(R\Gamma(X,\,\sE^r\x{\widetilde i}\sE^{r-1}\oplus\sE')\otimes_AA_v\Big)\simeq\det\Big(R\Gamma(X,\,\sE\x{i}\sE')\otimes_AA_v\Big);\\
&& \xi_i:\bigotimes_{w\in S}\det\Big(\sM^r_w\x{\widetilde i}\sM^{r-1}\oplus\sM'_w\Big)\simeq\bigotimes_{w\in S}\det\Big(\sM_w\x{i}\sM'_w\Big);\\
&&\eta_{i-j}:\det\Big(R\Gamma(X,\,\sE^r\x{\widetilde i-\widetilde j}\sE^{r-1}\oplus\sE')\otimes_AA_v\Big)\simeq\det\Big(R\Gamma(X,\,\sE\x{i-j}\sE')\otimes_AA_v\Big);\\
&& \xi_{i-j}:\bigotimes_{w\in S}\det\Big(\sM^r_w\x{\widetilde i-\widetilde j}\sM_w^{r-1}\oplus\sM'_w\Big)\simeq\bigotimes_{w\in S}\det\Big(\sM_w\x{i-j}\sM'_w\Big).
\end{eqnarray*}
By the functoriality of determinants, we have
\begin{eqnarray*}
\xi_i\widetilde\alpha=\alpha\eta_i\hbox{ and }\xi_{i-j}\widetilde\beta=\beta\eta_{i-j}.
\end{eqnarray*}
By Lemma \ref{d}, we have a diagram
\[\xymatrix{ \det\Big(R\Gamma(X,\,\sE^r\x{\widetilde i}\sE^{r-1}\oplus\sE')\Big)\ar[rr]^{\det(\sigma)\det({\rm id})^{-1}}\ar[d]^{\widetilde\varphi}
&&\det\Big(R\Gamma(X,\,\sE^r\x{{\rm id}\oplus i}\sE^{r-1}\oplus\sE')\Big)\ar[d]\ar[r]&\det\Big(R\Gamma(X,\,\sE\x{i}\sE')\Big)\ar[d]^\varphi\\
\det\Big(R\Gamma(X,\,\sE^r\x{\widetilde i-\widetilde j}\sE^{r-1}\oplus\sE')\Big)\ar[rr]^{\det(\phi_1)\det(\phi_2)^{-1}}&&
\det\Big(R\Gamma(X,\,\sE^r\x{{\rm id}\oplus i-j}\sE^{r-1}\oplus\sE')\Big)\ar[r]&\Big(R\Gamma(X,\,\sE\x{i-j}\sE')\Big).}\]
The right square commutes trivially. The left square commutes from the fact $\det(\sigma)=\det(\phi_1)=\det(\phi_2)=1$.
This proves  $\eta_{i-j}\widetilde\varphi=\varphi\eta_i$. Similarly, we have  $\xi_{i-j}\widetilde\psi=\psi\xi_i$.

By Theorem 5.1 of \cite{L}, we have
$$L_{v}^*(X-S,\,(\sE^r,\sE^{r-1}\oplus\sE',\widetilde i,\widetilde j))=[\widetilde\psi\widetilde\alpha:\widetilde \beta\widetilde\varphi]\in A_v^\times.$$
By (\ref{xiu}), we have
\begin{eqnarray*}
&&L_{v}^*(X-S,\,(\sE,\sE',i,j))\\
&=&L_{v}^*(X-S,\,(\sE^r,\sE^{r-1}\oplus\sE',\widetilde i,\widetilde j))\\
&=&[\widetilde\psi\widetilde\alpha:\widetilde \beta\widetilde\varphi]
=[\xi_{i-j}\widetilde\psi\widetilde\alpha:\xi_{i-j}\widetilde \beta\widetilde\varphi]
=[\psi\xi_i\widetilde\alpha:\beta\eta_{i-j}\widetilde\varphi]=[\psi\alpha\eta_i:\beta\varphi\eta_i]=[\psi\alpha:\beta\varphi].
\end{eqnarray*}
\end{proof}

\begin{lem}\label{j}
Keep the assumption of Lemma \ref{e}. Suppose for any $w\in S$, there exists $A_v$-linear isomorphisms
$\exp:\sM_w\simeq\sM_w$ and $\exp:\sM_w'\simeq\sM_w'$ which satisfy the following three conditions.
\begin{enumerate}
\item  $\exp\circ i=i-j\circ\exp:\sM_w\to\sM_w'$.
\item For any $t\in\N$, $\exp(w^t\sM_w)=w^t\sM_w$, $\exp(w^t\sM_w')=w^t\sM_w'$, $(\exp-{\rm id})(w^t\sM_w)\subset w^{t+1}\sM_w$ and $(\exp-{\rm id})(w^t\sM'_w)\subset w^{t+1}\sM'_w$.
\item For any $s\in\N$,  $(\exp-{\rm id})(w^t\sM_w)\subset w^{t+s}\sM_w$ and $(\exp-{\rm id})(w^t\sM'_w)\subset w^{t+s}\sM'_w$ for $t$ large enough.
\end{enumerate}
Let $\log:\sM_w\simeq\sM_w$ and $\log:\sM_w'\simeq\sM_w'$ be the inverse map of $\exp$. Let $\bigoplus_{w\in S}\sM_w'\x{\pi}C$ be the cokernel of $\bigoplus_{w\in S}\sM_w\x{i}\bigoplus_{w\in S}\sM_w'$. Denote by $\iota$ the natural map $R\Gamma(X,\,\sE')\otimes_AA_v\to\bigoplus_{w\in S}\sM_w'$.
Then we have a commutative diagram
\[\xymatrix{R\Gamma(X,\,\sE)\otimes_AA_v\ar[r]\ar@<0.5ex>[d]_{i-j~~}\ar@<-0.5ex>[d]
^{~~i}&\bigoplus_{w\in S}\sM_w\ar[r]\ar@<0.5ex>[d]_{i-j~~}\ar@<-0.5ex>[d]^{~~i}&{\rm Hom}_{A_v}(\sV,\,A_v)\ar@<0.5ex>[d]_{i-j~~}\ar@<-0.5ex>[d]^{~~i}\\
R\Gamma(X,\,\sE')\otimes_AA_v\ar[r]^\iota\ar@<0.5ex>[d]_{\pi\log\iota~~}\ar@<-0.5ex>[d]^{~~\pi\iota}&\bigoplus_{w\in S}\sM_w'\ar[r]\ar@<0.5ex>[d]_{\pi\log~~}\ar@<-0.5ex>[d]^{~~\pi}
&{\rm Hom}_{A_v}(\sV',\,A_v)\\
C\ar@{=}[r]&C},\]
whose four vertical triangles are distinguished.
They define four isomorphisms
\begin{eqnarray}
\label{1}&&\det(C)\simeq\det\Big(R\Gamma(X,\,\sE\x{i-j}\sE')\otimes_AA_v\Big)^{-1};\\
\label{2}&&\det(C)\simeq\det\Big(R\Gamma(X,\,\sE\x{i}\sE')\otimes_AA_v\Big)^{-1};\\
&&\det(C)\simeq\bigotimes_{w\in S}\det\Big(\sM_w\x{i-j}\sM_w'\Big)^{-1};\\
&&\det(C)\simeq\bigotimes_{w\in S}\det\Big(\sM_w\x{i}\sM_w'\Big)^{-1}.
\end{eqnarray}
Let $\delta$ be the composition
$$\det(C)\simeq\det\Big(R\Gamma(X,\,\sE\x{i}\sE')\otimes_AA_v\Big)^{-1}\x{\varphi^{-1}}\det\Big(R\Gamma(X,\,\sE\x{i-j}\sE')\otimes_AA_v\Big)^{-1}\simeq\det(C).$$
We have
$$L_{v}^*(X-S,\,(\sE,\sE',i,j))=[\delta:1]\in A_v^\times.$$
\end{lem}

\begin{proof}
Let $\mu$ be the composition
$$\det(C)\simeq\bigotimes_{w\in S}\det\Big(\sM_w\x{i}\sM_w'\Big)^{-1}\x{\psi^{-1}}\bigotimes_{w\in S}\det\Big(\sM_w\x{i-j}\sM_w'\Big)^{-1}\simeq\det(C).$$
By Lemma \ref{e}, we have $L_{v}^*(X-S,\,(\sE,\sE',i,j))=[\delta:\mu].$ To prove this lemma, it suffices to show that $\mu$ is the identity map. Recall the construction of $\psi$, we only need to prove for each $n\in\N$ the composition $\mu_n$
\begin{eqnarray*}\det(C/v^nC)&\simeq& \det\Big(\sM_w/w^t\sM_w+v^n\sM_w\x{i}(\sM_w'/v^n\sM_w')/i(w^t(\sM_w/v^n\sM_w))\Big)^{-1} \\
&\x{(\ref{c})^{-1}}&\det\Big(\sM_w/w^t\sM_w+v^n\sM_w\x{i-j}(\sM_w'/v^n\sM_w')/(i-j)(w^t(\sM_w/v^n\sM_w))\Big)^{-1} \\
&\simeq&\det(C/v^nC)
\end{eqnarray*}
is the identity map for $t$ large enough.
 By the commutative diagram
\[\xymatrix{0\ar[r]&\sM_w/w^t\sM_w+v^n\sM_w\ar[r]^i\ar[d]^\exp&\sM_w'/i(w^t\sM_w)+v^n\sM_w'\ar[rr]^\pi\ar[d]^\exp&&C/v^nC\ar@{=}[d]\ar[r]&0\\
0\ar[r]&\sM_w/w^t\sM_w+v^n\sM_w\ar[r]^{i-j}&\sM_w'/i(w^t\sM_w)+v^n\sM_w'\ar[rr]^{\pi\log}&&C/v^nC\ar[r]&0,}\]
for $t$ large enough, we have
$$[\mu_n:1]=\det\Big(\exp,\;\sM_w/w^t\sM_w+v^n\sM_w\Big)^{-1}\det\Big(\exp,\;\sM_w'/i(w^t\sM_w)+v^n\sM_w'\Big)\in A_v/v^nA_v.$$
For $0\leq s\leq t$, let $V_s=w^s(\sM_w/v^n\sM_w)/w^t(\sM_w/v^n\sM_w).$  Then we get a filtration $V_0\supset\cdots\supset V_t=0$ of $\sM_w/w^t\sM_w+v^n\sM_w$
by finite free $A_v/v^nA_v$-modules. By condition (2), $(\exp-{\rm id})(V_s)\subset V_{s+1}$ and hence $\det\Big(\exp,\;\sM_w/w^t\sM_w+v^n\sM_w\Big)=1$.
Similarly, $ \det\Big(\exp,\;\sM'_w/w^t\sM'_w+v^n\sM'_w\Big)=1$.
By Lemma 4.3 of \cite{L}, for each $n$, there exist $s\in\N$ such that $w^{t+s}(\sM_w'/v^n\sM_w')\subset i(w^t(\sM_w/v^n\sM_w))$ for $t$ large enough. By condition (3), we have
$$(\exp-{\rm id})(w^t(\sM_w'/v^n\sM_w'))\subset w^{t+s}(\sM_w'/v^n\sM_w')\subset i(w^t(\sM_w/v^n\sM_w)).$$
for $t$ large enough. Thus for $t\gg0$,
$$\det\Big(\exp,\;w^t(\sM_w'/v^n\sM_w')/i(w^t(\sM_w/v^n\sM_w))\Big)=1.$$
By the short exact sequence
$$0\to w^t(\sM_w'/v^n\sM_w')/i(w^t(\sM_w/v^n\sM_w))\to(\sM_w'/v^n\sM_w')/i(w^t(\sM_w/v^n\sM_w))\to\sM_w'/w^t\sM_w'+v^n\sM_w'\to0$$
of finite free $A_v/v^nA_v$-modules, we have
\begin{eqnarray*}
&&\det\Big(\exp,\;(\sM_w'/v^n\sM_w')/i(w^t(\sM_w/v^n\sM_w))\Big)\\
&=&\det\Big(\exp,\;w^t(\sM_w'/v^n\sM_w')/i(w^t(\sM_w/v^n\sM_w))\Big)
\det\Big(\exp,\;\sM_w'/w^t\sM_w'+v^n\sM_w'\Big)=1 \end{eqnarray*}
for $t\gg0$. This completes the proof of the lemma.
\end{proof}

\section{Proof of Theorem \ref{18} }
\noindent Before proving Theorem \ref{18}, we state three useful lemmas.
\begin{lem}\label{fang}
Let $S\to S'$ be a flat homomorphism of commutative rings.
Suppose $(C_1',C_1,C_1'',u_1,v_1,w_1)$ and $(C_2',C_2,C_2'',u_2,v_2,w_2)$ are distinguished triangles in $D^{\rm perf}(S)$ and $(C',C,C'',u,v,w)$ is a distinguished
triangle in $D^{\rm perf}(S')$. Any isomorphisms
$$\gamma':{\rm det}_S(C_1')\simeq{\rm det}_S(C_2'),\;\gamma:{\rm det}_S(C_1)\simeq{\rm det}_S(C_2)\hbox{ and }
\gamma'':{\rm det}_S(C_1'')\simeq{\rm det}_S(C_2'')$$
induces isomorphisms
\begin{eqnarray*}
&&\gamma'\otimes1:{\rm det}_{S'}(C_1'\otimes_SS')\simeq{\rm det}_{S'}(C_2'\otimes_SS'),\;
\gamma\otimes1:{\rm det}_{S'}(C_1\otimes_SS')\simeq{\rm det}_{S'}(C_2\otimes_SS'),\\
&&\gamma''\otimes1:{\rm det}_{S'}(C_1''\otimes_SS')\simeq{\rm det}_{S'}(C_2''\otimes_SS').
\end{eqnarray*}
For any isomorphism of triangles
\begin{eqnarray*}
(C_1'\otimes_SS',C_1\otimes_SS',C_1''\otimes_SS',u_1\otimes1,v_1\otimes1,w_1\otimes1)\x{(f_1',f_1,f_1'')}(C',C,C'',u,v,w);\\
(C_2'\otimes_SS',C_2\otimes_SS',C_2''\otimes_SS',u_2\otimes1,v_2\otimes1,w_2\otimes1)\x{(f_2',f_2,f_2'')}(C',C,C'',u,v,w),
\end{eqnarray*}
we have
$$[\det(f_2)(\gamma\otimes1)\det(f_1)^{-1}:1][1:\det(f_2')(\gamma'\otimes1)\det(f_1')^{-1}][1:\det(f_2'')(\gamma''\otimes1)\det(f_1'')^{-1}]\in S^\times\subset S'^\times.$$

Suppose furthermore $S$ and $S'$ are regular rings. Fix an isomorphism $\gamma_s:\det_S(H^sC_1)\simeq\det_S(H^sC_2)$ for each $s\in\Z$. Then we have
$$[\det(f_2)(\gamma\otimes1)\det(f_1)^{-1}:1]\prod_{s\in\Z}[\det(H^sf_2)(\gamma_s\otimes1)\det(H^sf_1)^{-1}:1]^{(-1)^{s+1}}\in S^\times.$$
\end{lem}

For any $\sO_X$-module $\sF$, let $\sF[t^{-1}]$, $\sF[t]$ and $\sF((t^{-1}))$ be the pull back of $\sF$ under the projections $X\otimes k[t^{-1}]\to X$
, $X\otimes k[t]\to X$ and $X\widehat\otimes k((t^{-1}))\to X$, respectively.
For any $k$-vector space $M$, let $M[t]=M\otimes_kk[t]$ and $M((t^{-1}))=M\widehat\otimes_kk((t^{-1})).$

\begin{lem}\label{k} Let $M$ be a $k$-vector space with a $k$-linear transformation $\theta$.
Then $\theta$ defines a $k[t]$-module structure on $M$.

(1) We have a short exact sequence
$$0\to M[t]\x{t-\theta}M[t]\x{p}M\to 0$$
of $k[t]$-modules, where $p(\sum_{s}m_st^s)=\sum_s\theta^s(m_s)$ for any $m_s\in M$.


(2) Consider the direct system
$(M_s,\,\theta)_{s\in\N}$ where $M_s=M$ and the morphisms $M_s\to M_{s+1}$ are all $\theta$. Let $p:M[t]=\bigoplus_{s\in\N}M_s\to \varinjlim (M_s,\theta)$ be
the natural surjective map.
We get a short exact sequence
$$0\to M[t]\x{1-t\theta}M[t]\x{p}\varinjlim_s (M_s,\theta)\to0$$
of $k[t]$-modules.
\end{lem}
\begin{proof}
(1) follows form Proposition 3, \cite{L1}.


(2) The injectivity follows form the fact $(1-t\theta)(m)-m\in t^{s+1}M[t]$ for any $m\in t^sM[t]$.
For any $0\neq m=\sum_{s=a}^bm_st^s$ with $m_a,\,m_b\neq0$, define the length $l(m)$ of $m$ to be $b-a$. If $m\in\ker(p)$, we prove that
$m\in(1-t\theta)M[t]$ by induction on $l(m)$. If $l(m)=0$, then $m=m_at^a$ for some $0\neq m_a\in M$ and $a\in\Z$. The condition $p(m_at^a)=0$
means that $\theta^c(m_a)=0$ for some $c\in\N$. Thus $m_at^a=(1-t\theta)\sum_{e=0}^c(t\theta)^e m_st^s.$ If $l(m)>0$, then
$l(m-(1-t\theta)m_at^a)<b-a$ and $m-(1-t\theta)m_at^a\in\ker(p)$. By induction hypothesis, $m-(1-t\theta)m_at^a\in(1-t\theta)M[t]$. So $m\in(1-t\theta)M[t]$.
\end{proof}

\begin{lem}\label{17} Let $M$ be a finite $k[t]$-module. Any isomorphism $\gamma:k[t]\simeq{\rm det}_{k[t]}(M)$ induces an isomorphism $\gamma\otimes1: k((t^{-1}))\simeq\det\Big(M\otimes_{k[t]}k((t^{-1}))\Big)$.
Since $M$ is finite, then $M\otimes_{k[t]}k((t^{-1}))=0$ and we have a canonical isomorphism $\rho:\det\Big(M\otimes_{k[t]}k((t^{-1}))\Big)\simeq k((t^{-1})).$
Then we have $[\rho(\gamma\otimes1):1]|M|\in k^\times\subset k((t^{-1}))^\times$.
\end{lem}

Recall that $R$ is the integral closure of of $k[t]$ in the function field $K$ of $X$ and $K_\infty=K\otimes_{k(t)}k((t^{-1}))$ and $Y={\rm Spec}\;R$.
Let $E$ be an abelian $t$-module over $R$ of dimension $n$ defined by
$$\varphi_E(t)=A_0+A_1\tau+\ldots A_r\tau^r\in M_n(R)\{\tau\}$$
such that $(A_0-zI_n)^n=0$ where $z$ is the image of $t$ in $R$.

Choose $e\in\N$ such that $A_0\in M_n(\Gamma(X,\sO_X(e\infty)))$. Then $e\geq1$. Let $\sO_\infty=\prod_{w\in X-Y}\sO_w$.
We get maps $1-t^{-1}A_0$ and $1-t^{-1}\sum_{s=0}^rA_s\tau^s:\sO_X(-d\infty)^n((t^{-1}))\to\sO_X((e-d)\infty)^n((t^{-1}))$ and $z^{-d}\sO_\infty^n((t^{-1}))\to z^{e-d}\sO_\infty^n((t^{-1}))$ for $d\gg0$.

Since $\lim_{s\to\infty}\frac{val(e_s)}{q^s}=+\infty$, we have $0\leq\sup_s-val(e_s)<+\infty$. Fix $c\in\N$ such that $c\geq\sup_{s\in\N}-val(e_s)+e$
and $A_1,\ldots,A_r\in M_n(\Gamma(X,\sO_X(c\infty)))$.
For any $x\in K_\infty^n$ with $val(x)>c-e$ and for any $s\geq1$, we have
$$val(e_sx^{(q^s)})=val(e_s)+q^sval(x)\geq val(e_s)+2val(x)>val(x).$$
This shows $val(\exp_E(x)-x)>val(x)$ for any $x\in K_\infty^n$ such that $val(x)>c-e$. Thus $\exp_E(z^{-d}\sO_\infty^n)\subset z^{-d}\sO_\infty^n$ and
the $\exp_E$ induces the identity map on $z^{-d}\sO_\infty^n/z^{-1-d}\sO_\infty^n$ for $d>c-e$. Thus $\exp_E$ defines an automorphism on $z^{-d}\sO_\infty^n$ for $d\gg0$.
It induces a $k((t^{-1}))$-linear automorphism on $z^{-d}\sO_\infty^n((t^{-1}))$ for $d\gg0$ which is also denoted by $\exp_E$. Then for $d\gg0$,
$$\exp_E\circ(1-t^{-1}A_0)=(1-t^{-1}\sum_{s=0}^rA_s\tau^s)\circ\exp_E:z^{-d}\sO_\infty^n((t^{-1}))\to z^{e-d}\sO_\infty^n((t^{-1})).$$
Let $\log_E$ be the inverse map of $\exp_E$. The above argument also shows that for any $s\in\N$, $(\exp_E-{\rm id})(z^{-d}\sO_\infty^n)\subset z^{-d-s}\sO_\infty^n$ for $d\gg0$.

Let $z^{e-d}\sO_\infty^n((t^{-1}))\x{\pi}C$ be the cokernel of $z^{-d}\sO_\infty^n((t^{-1}))\x{1-t^{-1}A_0}z^{e-d}\sO_\infty^n((t^{-1}))$.
Then $z^{e-d}\sO_\infty^n((t^{-1}))\x{\pi\log_E}C$ is the cokernel of $z^{-d}\sO_\infty^n((t^{-1}))\x{1-t^{-1}\sum_{s=0}^rA_s\tau^s}z^{e-d}\sO_\infty^n((t^{-1}))$.
Let $\iota:R\Gamma(X,\sO_X((e-d)\infty))^n((t^{-1}))\to z^{e-d}\sO_\infty^n((t^{-1}))$ be the natural map.
For any $d\in\Z$, we have a distinguished triangle
$$R\Gamma(X,\,\sO_X(-d\infty))\to z^{-d}\sO_\infty\to\frac{K_\infty}{R}.$$
Since $1-t^{-1}A_0$ and $1-t^{-1}\sum_{s=0}^rA_s\tau^s$ on $\big(\frac{K_\infty}{R}\big)^n((t^{-1}))$ are isomorphic, we have a commutative diagram
\[\xymatrix{R\Gamma(X,\,\sO_X(-d\infty))^n((t^{-1}))\ar[r]\ar@<0.5ex>[d]_{1-t^{-1}\sum_{s=0}^rA_s\tau^s~~}\ar@<-0.5ex>[d]
^{~~1-t^{-1}A_0}&z^{-d}\sO_\infty^n((t^{-1}))\ar[rr]\ar@<0.5ex>[d]_{1-t^{-1}\sum_{s=0}^rA_s\tau^s~~}\ar@<-0.5ex>[d]^{~~1-t^{-1}A_0}&&\big(\frac{K_\infty}{R}\big)^n((t^{-1}))
\ar@<0.5ex>[d]_{1-t^{-1}\sum_{s=0}^rA_s\tau^s~}\ar@<-0.5ex>[d]^{~1-t^{-1}A_0}\\
R\Gamma(X,\,\sO_X((e-d)\infty))^n((t^{-1}))\ar[r]\ar@<0.5ex>[d]_{\pi\log_E\iota~~}\ar@<-0.5ex>[d]^{~~\pi\iota}&z^{e-d}\sO_\infty^n((t^{-1}))
\ar[rr]\ar@<0.5ex>[d]_{\pi\log_E~~}\ar@<-0.5ex>[d]^{~~\pi}
&&\big(\frac{K_\infty}{R}\big)^n((t^{-1}))\\
C\ar@{=}[r]&C,}\]
whose four vertical triangles are distinguished.
\begin{lem}\label{jiang}
For $d\gg0$, we have $L(E/R)=[\delta_1:1]$ where $\delta_1$ is the composition
\begin{eqnarray*}
\det(C)&\simeq&\det\Big(R\Gamma\Big(X,\,\sO_X(-d\infty)^n((t^{-1}))\x{1-t^{-1}A_0}\sO_X((e-d)\infty)^n((t^{-1}))\Big)\Big)^{-1}\\
&\simeq&\det\Big(R\Gamma\Big(X,\,\sO_X(-d\infty)^n\Big)((t^{-1}))\Big)^{-1}\otimes\det\Big(R\Gamma\Big(X,\,\sO_X((e-d)\infty)^n\Big)((t^{-1}))\Big)\\
&\simeq&\det\Big(R\Gamma\Big(X,\,\sO_X(-d\infty)^n((t^{-1}))\x{1-t^{-1}\sum_{s=0}^rA_s\tau^s}\sO_X((e-d)\infty)^n((t^{-1}))\Big)\Big)^{-1}\simeq
\det(C).
\end{eqnarray*}
\end{lem}
\begin{proof}
For $d\gg0$, the diagram
$$\sO_X(-d\infty)^n[t^{-1}]\x{1-t^{-1}A_0}\sO_X((e-d)\infty)^n[t^{-1}]\xleftarrow{t^{-1}\sum_{s=1}^rA_s\tau^s}\sO_X(-d\infty)^n[t^{-1}]$$
satisfies the assumptions of Lemma \ref{j} at the place $t^{-1}$ of $k[t^{-1}]$. By Lemma \ref{j}, the $\infty$-adic special $L$-value
$$\prod_{\mathfrak m\in{\rm Max}(R)}{\rm det}_{k[[t^{-1}]]}\Big(1-(1-t^{-1}A_0)^{-1}\sum_{s=1}^rt^{-1}A_s\tau^s,\;(R/\mathfrak m)^n[[t^{-1}]]\Big)^{-1}$$ of the above diagram away from $X-Y$ converges to $[\delta_1:1]$. Then the lemma holds by
\begin{eqnarray*}
&&\prod_{\mathfrak m\in{\rm Max}(R)}{\rm det}_{k[[t^{-1}]]}\Big(1-(1-t^{-1}A_0)^{-1}\sum_{s=1}^rt^{-1}A_s\tau^s,\;(R/\mathfrak m)^n[[t^{-1}]]\Big)^{-1}\\
&=&\prod_{\mathfrak m\in{\rm Max}(R)}\frac{{\rm det}_{k[[t^{-1}]]}\Big(1-t^{-1}A_0,\;(R/\mathfrak m)^n[[t^{-1}]]\Big)}{{\rm det}_{k[[t^{-1}]]}\Big(1-t^{-1}\sum_{s=0}^rA_s\tau^s,\;(R/\mathfrak m)^n[[t^{-1}]]\Big)}\\
&=&\prod_{\mathfrak m\in{\rm Max}(R)}\frac{{\rm det}_{k[t]}\Big(t-A_0,\;(R/\mathfrak m)^n[t]\Big)}{{\rm det}_{k[t]}
\Big(t-\sum_{s=0}^rA_s\tau^s,\;(R/\mathfrak m)^n[t]\Big)}\\
&=&\prod_{\mathfrak m\in{\rm Max}(R)}\frac{|{\rm Lie}(E)(R/\mathfrak m)|}{|E(R/\mathfrak m)|}=L(E/R)\in 1+t^{-1}k[[t^{-1}]].
\end{eqnarray*}
\end{proof}

For $d\gg0$, we also have two distinguished triangles
\begin{eqnarray}\label{f1}
&&R\Gamma\Big(X,\,\sO_X(-d\infty)^n\Big)((t^{-1}))\x{t-A_0}R\Gamma\Big(X,\,\sO_X((e-d)\infty)^n\Big)((t^{-1}))\x{\pi\iota}C;\\\label{f2}
&&R\Gamma\Big(X,\,\sO_X(-d\infty)^n\Big)((t^{-1}))\x{t-\sum_{s=0}^rA_s\tau^s}R\Gamma\Big(X,\,\sO_X((e-d)\infty)^n\Big)((t^{-1}))\x{\pi\log_E\iota}C.
\end{eqnarray}
By Lemma \ref{jiang}, we have the following.
\begin{cor}\label{c1}
We have $L(E/R)=[\delta_2:1]$ where $\delta_2$ is the composition
\begin{eqnarray*}
\det(C)\simeq\det\Big(R\Gamma\Big(X,\,\sO_X(-d\infty)^n\Big)((t^{-1}))\Big)^{-1}\bigotimes\det\Big(R\Gamma\Big(X,\,\sO_X((e-d)\infty)^n\Big)((t^{-1}))\Big)\simeq\det(C)
\end{eqnarray*}
where the two isomorphisms are defined by (\ref{f1}) and (\ref{f2}), respectively.
\end{cor}

By Lemma \ref{fu}, we can define a $k((t^{-1}))$-linear map $q$
$$q:z^{e-d}\sO_\infty^n((t^{-1}))\to{\rm Lie}(E)(K_\infty),\;\;\sum_{s\ll+\infty}x_st^s\mapsto\sum_{s\ll+\infty}A_0^s(x_s)$$
for any $x_s\in z^{e-d}\sO_\infty^n$. For any $x\in{\rm Lie}(E)(K_\infty)$, $z^{-q^{ns}}x\in z^{e-d}\sO_\infty^n$ for $s\gg0$. Then
$q(z^{-q^{ns}}xt^{q^{ns}})=A_0^{q^{ns}}(z^{-q^{ns}}x)=x$. This shows $q(z^{e-d}\sO_\infty^n[t])={\rm Lie}(E)(K_\infty)$.
The surjective map $q$ factors as $z^{e-d}\sO_\infty^n((t^{-1}))\x{\pi}C\x{p}{\rm Lie}(E)(K_\infty)$.
Let $z^{e-d}\sO_\infty^n[t]\x{\pi_1}D$ be the cokernel of $z^{-d}\sO_\infty^n[t]\x{t-A_0}z^{e-d}\sO_\infty^n[t].$
Then $z^{e-d}\sO_\infty^n[t]\x{\pi_1\log_E}D$ is the cokernel of $z^{-d}\sO_\infty^n[t]\x{t-\sum_{s=0}^rA_s\tau^s}z^{e-d}\sO_\infty^n[t].$
There exists a unique $k[t]$-linear map $\eta:D\to C$ such that the following diagram commutes:
\[\xymatrix{z^{-d}\sO_\infty^n[t]\ar[r]^{t-A_0}\ar[d]&z^{e-d}\sO_\infty^n[t]\ar[rr]^{\pi_1}\ar[d]&&D\ar[d]^\eta\ar[rr]^{p\eta}&&{\rm Lie}(E)(K_\infty)\ar@{=}[d]\\
z^{-d}\sO_\infty^n((t^{-1}))\ar[r]^{t-A_0}&z^{e-d}\sO_\infty^n((t^{-1}))\ar[rr]^\pi &&C\ar[rr]^p&&{\rm Lie}(E)(K_\infty).}\]
Since $p\eta\pi_1$ is surjective, then so is $p\eta$ and $p$. By Lemma \ref{k} (1), we have two short exact sequences
\begin{eqnarray*}\label{10}
&&0\to\Big(\frac{K_\infty}{R}\Big)^n[t]\x{t-A_0}\Big(\frac{K_\infty}{R}\Big)^n[t]\x{}{\rm Lie}(E)\Big(\frac{K_\infty}{R}\Big)\to0;\\\label{11}
&&0\to\Big(\frac{K_\infty}{R}\Big)^n[t]\x{t-\sum_{s=0}^rA_s\tau^s}\Big(\frac{K_\infty}{R}\Big)^n[t]\x{}E\Big(\frac{K_\infty}{R}\Big)\to0
\end{eqnarray*}
Let $\bar p$ be the natural map ${\rm Lie}(E)(K_\infty)\to{\rm Lie}(E)(\frac{K_\infty}{R})$ and let $\overline\exp_E:{\rm Lie}(E)(K_\infty)\to E(\frac{K_\infty}{R})$ be the map induced by $\exp_E:{\rm Lie}(E)(K_\infty)\to E(K_\infty)$.
By the commutative diagram
\[\xymatrix{R\Gamma\big(X,\,\sO_X(-d\infty)^n\big)[t]\ar[r]\ar@<0.5ex>[d]_{t-\sum_{s=0}^rA_s\tau^s~~}\ar@<-0.5ex>[d]
^{~~t-A_0}&z^{-d}\sO_\infty^n[t]\ar[rr]\ar@<0.5ex>[d]_{t-\sum_{s=0}^rA_s\tau^s~~}\ar@<-0.5ex>[d]^{~~t-A_0}
&&\big(\frac{K_\infty}{R}\big)^n[t]\ar@<0.5ex>[d]_{t-\sum_{s=0}^rA_s\tau^s~~}\ar@<-0.5ex>[d]^{~~t-A_0}\\
R\Gamma\big(X,\,\sO_X((e-d)\infty)^n\big)[t]\ar[r]&z^{e-d}\sO_\infty^n[t]
\ar[rr]&&\big(\frac{K_\infty}{R}\big)^n[t],}\]
we get two quasi-isomorphisms
\begin{eqnarray}\label{f3}
&&R\Gamma\Big(X,\,\sO_X(-d\infty)^n[t]\x{t-A_0}\sO_X((e-d)\infty)^n[t]\Big)[1]
\simeq\Big(D\x{\bar pp\eta}{\rm Lie}(E)\Big(\frac{K_\infty}{R}\Big)\Big);\\\label{f4}
&&R\Gamma\Big(X,\,\sO_X(-d\infty)^n[t]\x{t-\sum_{s=0}^rA_s\tau^s}\sO_X((e-d)\infty)^n[t]\Big)[1]
\simeq\Big(D\x{\overline{\exp}_Ep\eta}E\Big(\frac{K_\infty}{R}\Big)\Big).
\end{eqnarray}
Since $X$ is projective, then (\ref{f3}) and (\ref{f4}) show that
\begin{eqnarray}\label{li}\Big(D\x{\bar pp\eta}{\rm Lie}(E)\Big(\frac{K_\infty}{R}\Big)\Big)\hbox{ and }\Big(D\x{\overline{\exp}_Ep\eta}E\Big(\frac{K_\infty}{R}\Big)\Big)\in D^{\rm perf}(k[t]).\end{eqnarray}
By (\ref{f3}) and (\ref{f4}), the composition
\begin{eqnarray*}
&&\det\Big(R\Gamma\Big(X,\,\sO_X(-d\infty)^n[t]\x{t-A_0}\sO_X((e-d)\infty)^n[t]\Big)\Big)^{-1}\\
&\simeq&\det\Big(R\Gamma\Big(X,\,\sO_X(-d\infty)^n[t]\Big)\Big)^{-1}\bigotimes\det\Big(R\Gamma\Big(X,\,\sO_X((e-d)\infty)^n[t]\Big)\Big)\\
&\simeq&\det\Big(R\Gamma\Big(X,\,\sO_X(-d\infty)^n[t]\x{t-\sum_{s=0}^rA_s\tau^s}\sO_X((e-d)\infty)^n[t]\Big)\Big)^{-1}
\end{eqnarray*}
defines an isomorphism
\begin{eqnarray*}\label{f7}
\gamma:\det\Big(D\x{\bar pp\eta}{\rm Lie}(E)\Big(\frac{K_\infty}{R}\Big)\Big)\simeq\det\Big(D\x{\overline{\exp}_Ep\eta}E\Big(\frac{K_\infty}{R}\Big)\Big).
\end{eqnarray*}

By (\ref{f1}), (\ref{f3}), (\ref{f2}) and (\ref{f4}), we have two natural quasi-isomorphisms
\begin{eqnarray*}\label{f5}
&&g_1:\Big(D\x{\bar pp\eta}{\rm Lie}(E)\Big(\frac{K_\infty}{R}\Big)\Big)\otimes_{k[t]}k((t^{-1}))\simeq C;\\\label{f6}
&&g_2:\Big(D\x{\overline{\exp}_Ep\eta}E\Big(\frac{K_\infty}{R}\Big)\Big)\otimes_{k[t]}k((t^{-1}))\simeq C.
\end{eqnarray*}
By Corollary \ref{c1}, we have
\begin{eqnarray}\label{c2}
L(E/R)=[\det(g_2)\circ(\gamma\otimes 1)\circ\det(g_1)^{-1}:1]\in k((t^{-1}))^\times.
\end{eqnarray}

Let $\sK$ be the kernel of the surjective map $D\x{p\eta}{\rm Lie}(E)(K_\infty)$. We get two distinguished triangles
\begin{eqnarray}\label{f8}
&& \sK\to\Big(D\x{\bar pp\eta}{\rm Lie}(E)\Big(\frac{K_\infty}{R}\Big)\Big)\to\Big({\rm Lie}(E)(K_\infty)\x{\bar p}{\rm Lie}(E)\Big(\frac{K_\infty}{R}\Big)\Big);\\\label{f9}
&& \sK\to\Big(D\x{\overline{\exp}_Ep\eta}E\Big(\frac{K_\infty}{R}\Big)\Big)\to\Big({\rm Lie}(E)(K_\infty)\x{\overline{\exp}_E}E\Big(\frac{K_\infty}{R}\Big)\Big),
\end{eqnarray}
in $D(k[t])$. Let $L$ be the kernel of the surjective map $C\x{p}{\rm Lie}(E)(K_\infty)$.
We get a short exact sequence
\begin{eqnarray}\label{f10}
0\to L\to C\x{p}{\rm Lie}(E)(K_\infty)\to0
\end{eqnarray}
of $k((t^{-1}))$-vector spaces.
The map $\eta:D\to C$ induces a $k[t]$-linear map $\sK\to L$ and hence a $k((t^{-1}))$-linear map $\sK\otimes_{k[t]}k((t^{-1}))\x{f_1}L$.
The identity map of ${\rm Lie}(E)(K_\infty)$ as $k[t]$-modules induces a $k((t^{-1}))$-map ${\rm Lie}(E)(K_\infty)\otimes_{k[t]}k((t^{-1}))\to{\rm Lie}(E)(K_\infty)$.
Let $h_1$ and $h_2$ be the composition
\begin{eqnarray*}
&&h_1:\Big({\rm Lie}(E)(K_\infty)\x{\bar p}{\rm Lie}(E)\Big(\frac{K_\infty}{R}\Big)\Big)\otimes_{k[t]}k((t^{-1}))\to{\rm Lie}(E)(K_\infty)\otimes_{k[t]}k((t^{-1}))\to{\rm Lie}(E)(K_\infty);\\
&&h_2:\Big({\rm Lie}(E)(K_\infty)\x{\overline{\exp}_E}E\Big(\frac{K_\infty}{R}\Big)\Big)\otimes_{k[t]}k((t^{-1}))\to{\rm Lie}(E)(K_\infty)\otimes_{k[t]}k((t^{-1}))\to{\rm Lie}(E)(K_\infty).
\end{eqnarray*}
Since $A_0^{q^n}=z^{q^n}I_n$, we have ${\rm Lie}(E)(R)=R^n$ as $k[t^{q^n}]$-modules. Then ${\rm Lie}(E)(R)$ is a finite generated $k[t^{q^n}]$-module and hence it is a finite generated $k[t]$-module. Since $A_0$ is invertible on $K_\infty^n$, then ${\rm Lie}(E)(R)\subset{\rm Lie}(E)(K_\infty)$ is a torsion free $k[t]$-module.
The inclusion ${\rm Lie}(E)(R)\subset{\rm Lie}(E)(K_\infty)$ defines a $k((t^{-1}))$-linear map ${\rm Lie}(E)(R)\otimes_{k[t]}k((t^{-1}))\to{\rm Lie}(E)(K_\infty)$.
As $k((t^{-q^n}))$-vector spaces, we have
\begin{eqnarray*}
{\rm Lie}(E)(R)\otimes_{k[t]}k((t^{-1}))\simeq{\rm Lie}(E)(R)\otimes_{k[t^{q^n}]}k((t^{-q^n}))\simeq R^n\otimes_{k[t^{q^n}]}k((t^{-q^n}))\simeq K_\infty^n\simeq{\rm Lie}(E)(K_\infty).
\end{eqnarray*}
Then the $k((t^{-1}))$-linear map ${\rm Lie}(E)(R)\otimes_{k[t]}k((t^{-1}))\to{\rm Lie}(E)(K_\infty)$ is an isomorphism. This proves ${\rm Lie}(E)(R)$ is a lattice of ${\rm Lie}(E)(K_\infty)$.
Then $h_1$ is a quasi-isomorphism by the fact
$$\Big({\rm Lie}(E)(K_\infty)\x{\bar p}{\rm Lie}(E)\Big(\frac{K_\infty}{R}\Big)\Big)\simeq{\rm Lie}(E)(R)\in D^{\rm perf}(k[t]).$$
By (\ref{li}) and (\ref{f8}), $\sK\in D^{\rm perf}(k[t])$. By (\ref{li}) and (\ref{f9}),
$$\Big({\rm Lie}(E)(K_\infty)\x{\overline{\exp}_E}E\Big(\frac{K_\infty}{R}\Big)\Big)\in D^{\rm perf}(k[t]).$$
Consider the natural morphism of triangles
\begin{eqnarray*}
&&(\ref{f8})\otimes_{k[t]}k((t^{-1}))\x{(f_1,g_1,h_1)}(\ref{f10});\\\nonumber
&&(\ref{f9})\otimes_{k[t]}k((t^{-1}))\x{(f_1,g_2,h_2)}(\ref{f10}).
\end{eqnarray*}
Since $g_1$ and $h_1$ are isomorphisms, so is $f_1$. Since $g_2$ is an isomorphism, so is $h_2$.
The distinguished triangles (\ref{f8}) and (\ref{f9}) defines two isomorphisms
\begin{eqnarray*}
&& \det\Big({\rm Lie}(E)(K_\infty)\x{\bar p}{\rm Lie}(E)\Big(\frac{K_\infty}{R}\Big)\Big)\simeq\det(\sK)^{-1}\bigotimes\det\Big(D\x{\bar pp\eta}{\rm Lie}(E)\Big(\frac{K_\infty}{R}\Big)\Big);\\
&& \det\Big({\rm Lie}(E)(K_\infty)\x{\overline{\exp}_E}E\Big(\frac{K_\infty}{R}\Big)\Big)\simeq\det(\sK)^{-1}\bigotimes\det
\Big(D\x{\overline{\exp}_Ep\eta}E\Big(\frac{K_\infty}{R}\Big)\Big).
\end{eqnarray*}
The isomorphism $\gamma$ defines an isomorphism
\begin{eqnarray*}
\gamma_1:\det\Big({\rm Lie}(E)(K_\infty)\x{\bar p}{\rm Lie}(E)\Big(\frac{K_\infty}{R}\Big)\Big)\simeq\det\Big({\rm Lie}(E)(K_\infty)\x{\overline{\exp}_E}E\Big(\frac{K_\infty}{R}\Big)\Big).
\end{eqnarray*}
Applying Lemma \ref{fang} to $k[t]\hookrightarrow k((t^{-1}))$, we have
\begin{eqnarray}\label{xiao}
&&[\det(g_2)(\gamma\otimes1)\det(g_1)^{-1}:1]\\\nonumber
&=&[\det(f_1)({\rm id}_{\det(\sK)}\otimes1)\det(f_1)^{-1}:1][\det(h_2)(\gamma_1\otimes1)\det(h_1)^{-1}:1]\\\nonumber
&=&[\det(h_2)(\gamma_1\otimes1)\det(h_1)^{-1}:1]
\end{eqnarray}
We have
\begin{eqnarray*}\label{si}
&&H^0\Big({\rm Lie}(E)(K_\infty)\x{\bar p}{\rm Lie}(E)\Big(\frac{K_\infty}{R}\Big)\Big)={\rm Lie}(E)(R);\\\nonumber
&&H^0\Big({\rm Lie}(E)(K_\infty)\x{\overline{\exp}_E}E\Big(\frac{K_\infty}{R}\Big)\Big)=\exp_E^{-1}(E(R));\\\nonumber
&&H^1\Big({\rm Lie}(E)(K_\infty)\x{\overline{\exp}_E}E\Big(\frac{K_\infty}{R}\Big)\Big)=H(E/R);\\\nonumber
&&H^s\Big({\rm Lie}(E)(K_\infty)\x{\bar p}{\rm Lie}(E)\Big(\frac{K_\infty}{R}\Big)\Big)=H^s\Big({\rm Lie}(E)(K_\infty)\x{\overline{\exp}_E}E\Big(\frac{K_\infty}{R}\Big)\Big)=0,\hbox{ otherwise }.
\end{eqnarray*}
Since $h_2$ is an isomorphism and $\Big({\rm Lie}(E)(K_\infty)\x{\overline{\exp}_E}E\Big(\frac{K_\infty}{R}\Big)\Big)$ is perfect, then $\exp_E^{-1}(E(R))$ is a lattice of ${\rm Lie}(E)(K_\infty)$ and $H(E/R)$ is a finite $k[t]$-module.

Fix an isomorphism $\gamma_2:\det\Big({\rm Lie}(E)(R)\Big)\simeq\det\Big(\exp_E^{-1}(E(R))\Big)$ and $\gamma_3:k[t]\simeq\det\Big( H(E/R)\Big)$.
By Lemma \ref{fang}, Lemma \ref{17} and Definition \ref{16}, we have
\begin{eqnarray}\label{sun}
&&[\det(h_2)(\gamma_1\otimes1)\det(h_1)^{-1}:1]\\\nonumber
&=&[\det\big(H^0(h_2)\big)(\gamma_2\otimes1)\det\big(H^0(h_1)\big)^{-1}:1][\det\big(H^1(h_2)\big)(\gamma_3\otimes1)\det\big(H^1(h_1)\big)^{-1}:1]^{-1}\\\nonumber
&=&[{\rm Lie}(E)(R):\exp_E^{-1}(E(R))]\cdot|H(E/R)|\in k((t^{-1}))^\times/k^\times.
\end{eqnarray}
Since $L(E/R)$, $[{\rm Lie}(E)(R):\exp_E^{-1}(E(R))]$ and $|H(E/R)|$ are monic, then by (\ref{c2}), (\ref{xiao}), (\ref{sun}) and Definition \ref{16}, we have
$$L(E/R)=[{\rm Lie}(E)(R):\exp_E^{-1}(E(R))]\cdot|H(E/R)|\in k((t^{-1}))^\times.$$

\section{Proof of Theorem \ref{yang}}
Recall the ${\rm Spec}\,k[t]$-shtuka $\widetilde\sE$ on $X$, the complex $\sE^\bullet$ on $X\otimes k[t]$ and the sheaf $\mathfrak E$ on $X$ in Definition \ref{ke}.

Let $j_Y:Y\to X$ be the inclusion. Since $K_\infty=\prod_{w\in X-Y}K_w$, then $K_\infty^n$ can be viewed as a sheaf on $X$ supported on $X-Y$.
So is $z^{-d}\sO_\infty^n[t]$ and $z^{e-d}\sO_\infty^n[t]$.
Consider the morphism
\[\xymatrix{0\ar[r]&\sO_X(-d\infty)^n[t]\ar[r]\ar[d]^{t-\sum_{s=0}^rA_s\tau^s~~}&j_{Y*}\sO_Y^n[t]\times z^{-d}\sO_\infty^n[t]\ar[rr]^f\ar[d]^{t-\sum_{s=0}^rA_s\tau^s~~}&& K_\infty^n[t]\ar[d]^{t-\sum_{s=0}^rA_s\tau^s}\ar[r]&0
\\
0\ar[r]&\sO_X((e-d)\infty)^n[t]\ar[r]&j_{Y*}\sO_Y^n[t]\times z^{e-d}\sO_\infty^n[t]\ar[rr]^f&&K_\infty^n[t]\ar[r]&0,}\]
of short exact sequences, where $f(x,y)=x-y$.
Recall that $z^{e-d}\sO_\infty^n[t]\x{\pi_1\log_E}D$ is the cokernel of $z^{-d}\sO_\infty^n[t]\x{t-\sum_{s=0}^rA_s\tau^s}z^{e-d}\sO_\infty^n[t]$. By Lemma \ref{k} (1), $\sE^\bullet[1]$ is a sheaf and
\begin{eqnarray}\label{mi}
0\to\sE^\bullet[1]\to j_{Y*}E(\sO_Y)\times D\to E(K_\infty)\to0
\end{eqnarray}
is exact, where the last map is $(x,y)\mapsto x-\exp_Ep\eta(y)$.
Recall the definition of $\mathfrak E$, we get a short exact sequence
\begin{eqnarray}\label{ng}
0\to\mathfrak E\to j_{Y*}E(\sO_Y)\times {\rm Lie}(E)(K_\infty)\to E(K_\infty)\to 0,
\end{eqnarray}
where the last map is $(x,y)\mapsto x-\exp_E(y)$. The kernel of the surjective map $p\eta:D\to{\rm Lie}(E)(K_\infty)$ is $\sK$. By (\ref{mi}) and (\ref{ng}), $p\eta$ induces a surjective morphism
$\sE^\bullet[1]\to\mathfrak E$ with kernel $\sK$.
We have a long exact sequence
\begin{eqnarray}\label{lei}
o\to\sK\to H^1(X,\,\sE^\bullet)\to H^0(X,\,\mathfrak E)\to0\to H^2(X,\,\sE^\bullet)\to H^1(X,\,\mathfrak E)\to0.
\end{eqnarray}
By (\ref{f4}) and (\ref{f9}), we have a quasi-isomorphism
\begin{eqnarray*}
R\Gamma(X,\,\mathfrak E)\simeq\Big({\rm Lie}(E)(K_\infty)\x{\overline{\exp}_Ep\eta}E\Big(\frac{K_\infty}{R}\Big)\Big).
\end{eqnarray*}
This proves
\begin{eqnarray}\label{yu}
&&H^0(X,\,\mathfrak E)= \exp_E^{-1}(E(R));\\\nonumber
&&H^1(X,\,\mathfrak E)= H(E/R).
\end{eqnarray}
By Proposition 5 of \cite{L1},
we have
$${\rm Ext}(\textbf{1}_{X\otimes k[t]},\;\widetilde\sE)\simeq R\Gamma(X,\,\sE^\bullet).$$
By (\ref{lei}) and (\ref{yu}), we have a short exact sequence
$$0\to\sK\to{\rm Ext}^1(\textbf{1}_{X\otimes k[t]},\;\widetilde\sE)\to\exp_E^{-1}(E(R))\to0$$
and
\begin{eqnarray*}
&&{\rm Ext}^2(\textbf{1}_{X\otimes k[t]},\;\widetilde\sE)\simeq H(E/R);\\
&&{\rm Ext}^s(\textbf{1}_{X\otimes k[t]},\;\widetilde\sE)=0\hbox{ for }
s\neq 1\hbox{ and }2.
\end{eqnarray*}
Suppose $A_0\in M_n(\Gamma(X,\,\sO_X(\infty)))$ and let $e=1$. By $z^{-1}\Gamma(X,\,\sO_X(\infty))\subset\sO_\infty$, we have $z^{-1}A_0=I_n+Q$ for some $Q\in M_n(\sO_\infty)$ such that $Q^n=0$. Hence
$z^{-1}A_0$ and $zA_0^{-1}$ are invertible matrixes in $M_n(\sO_\infty)$. So $A_0^{-1}\in M_n(\sO_\infty)$. So $A_0$ defines an
isomorphism $z^{-d}\sO_\infty^n\simeq z^{1-d}\sO_\infty^n$ for each $d$.
Then
$$\varinjlim(z^{1-d}\sO_\infty^n,\,A_0^{-1})=\varinjlim(z^{1-d}\sO_\infty^n,\,A_0^{-qn})=\varinjlim(z^{1-d}\sO_\infty^n,\,z^{-qn})={\rm Lie}(E)(K_\infty).$$
Recall that $q$ is the map
$$z^{1-d}\sO_\infty^n[t]\to {\rm Lie}(E)(K_\infty),\;\;\sum_sx_st^s\mapsto \sum_sA_0^s(x_s)$$
for any $x_s\in z^{1-d}\sO_\infty^n$. By Lemma \ref{j} (3), we have a short exact sequence of $k[t]$-modules
$$0\to z^{1-d}\sO_\infty^n[t]\x{1-tA_0^{-1}}z^{1-d}\sO_\infty^n[t]\x{q}{\rm Lie}(E)(K_\infty)\to0.$$
By the isomorphism $A_0:z^{-d}\sO_\infty^n[t]\simeq z^{1-d}\sO_\infty^n[t]$, we have a short exact sequence of $k[t]$-modules
\begin{eqnarray}\label{8}
0\to z^{-d}\sO_\infty^n[t]\x{t-A_0}z^{1-d}\sO_\infty^n[t]\x{q}{\rm Lie}(E)(K_\infty)\to0.
\end{eqnarray}
Then $p\eta:D\to{\rm Lie}(E)(K_\infty)$ is an isomorphism and $\sK=0$. So
$${\rm Ext}^1(\textbf{1}_{X\otimes k[t]},\;\widetilde\sE)\simeq\exp_E^{-1}(E(R)).$$
This completes the proof of Theorem \ref{yang}.

\bibliographystyle{plain}

\end{document}